\documentclass[12pt]{article}
\usepackage{amsmath}
\usepackage{amssymb}
\usepackage{verbatim}
\newtheorem{thm}{Theorem}
\newtheorem{lemma}[thm]{Lemma}

\newtheorem{prop}[thm]{Proposition}
\newtheorem{definition}[thm]{Definition}

\newcommand{\dif}[0]{\ensuremath{\,\mathrm{d}}}

\newcommand{\e}{\varepsilon}
\renewcommand{\u}{u_{\varepsilon}}
\renewcommand{\v}{v_{\varepsilon}}
\newcommand{\V}{V_{\varepsilon}}

\DeclareMathOperator*{\dive}{div}

\begin{document}
\title{ THE TIME DERIVATIVE IN A SINGULAR PARABOLIC EQUATION}
\author{Peter Lindqvist}
\date{\footnotesize{Department of Mathematical Sciences\\Norwegian University of Science and Technology\\NO--7491 Trondheim, Norway}}
\maketitle

\begin{center}
  \texttt{Dedicated to Olli Martio on his seventy-fifth birthday}
\end{center}
 \bigskip
{\small \textsc{Abstract:}}\footnote{AMS classification 35K67,
35K92, 35B45} \textsf{We study the Evolutionary $p$-Laplace Equation in the singular case $1<p<2$. We prove that a weak solution has a time derivative (in Sobolev's sense) which is a function belonging (locally) to a $L^q$-space.}

\section{Introduction}
The regularity theory for parabolic partial differential equations of the type
$$\frac{\partial u}{\partial t} \,=\,\dive \mathbf{A}(x,t,u,\nabla u)$$
aims at establishing boundedness and continuity of the solution $u =u(x,t)$ and its gradient $$\nabla u = \Bigl(\dfrac{\partial u}{\partial x_{1}},\dots,\frac{\partial u}{\partial x_{n}}\Bigr).$$
The celebrated methods of DeGiorgi, Nash, and Moser do not directly treat the time derivative $u_t,$ which is regarded as merely a distribution. Yet, for many specific equations the time derivative is more than that, it is a  function in Lebesgue's theory. We shall prove that the solutions of the \emph{Evolutionary} $p$-\emph{Laplace Equation}
\begin{equation}
  \frac{\partial u}{\partial t} \,=\,\dive\left(|\nabla u|^{p-2}\nabla u \right)
\end{equation}
have a first order time derivative $u_t$  in Sobolev's sense. Thus the time derivative exists as a measurable function satisfying the definition
$$ \int\limits_0^T\!\!\!\int\limits_{\Omega}\u\nabla\phi\,dxdt\,=\, - \int\limits_0^T\!\!\!\int\limits_{\Omega}u_t\phi\,dxdt$$
for all test functions $\phi \in C^{\infty}_0(\Omega_T).$ Here $\Omega_T = \Omega \times (0,T)$ and $\Omega$ is a domain in the $n$-dimensional space $\mathbb{R}^n.$ The so-called \textsf{degenerate case} $p\geq 2$ (or slow diffusion case) was treated in [L1] and [L2] and now we shall focus our attention on the so-called \textsf{singular case} (or fast diffusion case) $1<p<2,$ which is much more demanding, because the operator
$$ \dive(|\nabla u|^{p-2}\nabla u) \,=\, |\nabla u|^{p-2}\Delta u + (p-2)|\nabla u|^{p-4} \sum_{i,j}^{n}\frac{\partial u}{\partial x_{i}}\frac{\partial u}{\partial x_{j}}\frac{\partial^2 u}{\partial x_{i}\partial x_{j}}$$
is undefined at the critical points $\nabla u = 0$ when $p<2.$ (It is known that the second derivatives $u_{x_ix_j}$ exist in the singular case, but the negative power $p-2$ spoils the formula.) --- We refer to the books [DB], [WYZ] about the Evolutionary $p$-Laplace Equation Equation.

Our method is to differentiate the regularized equation
\begin{equation}\label{regu}
  \frac{\partial u}{\partial t}\,=\, \dive\Bigl(\!\bigl\{|\nabla u|^2 + \e^2\bigr\}^{\frac{p-2}{2}}\nabla u\!\Bigr)
  \end{equation}
with respect to the $x$-variables and then to derive careful estimates which are passed over to the limit as $\e \to 0.$ The appearing identities are, of course, not new. The main formula to start from has been used for other purposes in [Y] and [WZY]. The case $p \geq \frac{3}{2}$ can be extracted from [Y]. See also [AMS] for systems.  Unfortunately, there is an extra complication when $p$ is small; in our proofs it appears in the range $1<p<\frac{3}{2}.$ To wit, the natural  definition that weak solutions belong only to
$$u\in C\bigl(0,T;L^2(\Omega)\bigr)\,\cap\, L^p\bigl(0,T;W^{1,p}(\Omega)\bigr)$$
is problematic, because it allows for \emph{un}bounded solutions, when $1<p<\frac{2n}{n+2}.$ See [DB] and [DH] for this striking phenomenon. Various attempts to deal with this situation are suggested in [BIV], where it is even proposed to assume that $u_t$ and $\dive(|\nabla u|^{p-2}\nabla u)$ both belong to the space $L^1.$ (Lemma III.3.6 in [DH] states that $u_t\in L^1_{loc}$ under the condition that  the B\'{e}nilan-Crandall estimate
$$u_t\,\leq \,\frac{1}{2-p}\,\frac{u}{t},$$
is valid, which requires some restrictions.)
The common definition is to add the condition 
$$u \in L^{r}(\Omega_T),\qquad \text{where}\quad p(n+r) > 2n$$
for some exponent $r$ in this range.
This extra assumption has the effect that the weak solutions become locally bounded. See [DH], in particular III.6 and III.7 for further information about this sharp condition.

We shall directly assume  that $\|u\|_{\infty} < \infty$, if $ 1<p<\frac{3}{2}.$ However, in the one-dimensional case ($n=1$) we have, without any extra hypothesis, a short proof that the time derivative is square summable.\footnote{Neither is any extra condition needed for $n=2$, since $p(\mathbf{2}+ r) > 2\cdot\mathbf{2},$ if $r=2.$}

\bigskip
\begin{thm} Suppose that $u$ is a weak solution of the equation $$u_t=\dive(|\nabla u|^{p-2}\nabla u)$$ in the domain $\Omega_T$. In the case $1<p<\frac{3}{2},\,\, n\geq 2,$ we make the \emph{extra} assumption that $u$ is bounded.  Then the time derivative $u_t$ exists in Sobolev's sense and $u_t \in L_{loc}^{\theta}(\Omega_T)$ for some $\theta > 1.$

  If $p \geq \frac{3}{2}$ or $n = 1$, we can take $\theta = 2.$

  If $1<p<\frac{3}{2}$ and $ n \geq 2,$ we have the restriction $1< \theta < \frac{1}{2-p}.$
\end{thm}

\bigskip
\noindent A quantitative proof is the object of this work. It is noteworthy that the proper regularity theory is not invoked.

\section{Preliminaries}

We use standard notation. See [DB] about  time dependent Sobolev Spaces.
Suppose that $\Omega$ is a bounded domain in $\mathbb{R}^n$ and consider the space$\times$time cylinder $\Omega_T = \Omega\times (0,T).$ We shall always keep $1 < p \leq 2,$ although many formulas are valid also for $p > 2.$ Denote $|\mathbb{D}^2u|^2 = \sum u^2_{x_{i}x_{j}}.$ Once and for all, we fix a test function $\zeta \in C_0^{\infty}(\Omega),\,\,0\leq\zeta\leq1.$ In the sequel, the constants in the estimates can depend on $\|\zeta_t\|_{\infty}$ and $\|\nabla \zeta\|_{\infty}.$
\bigskip
\begin{definition} Assume that  $u \in C\bigl(0,T;L^2(\Omega)\bigr) \cap L^p\bigl(0,T;W^{1,p}(\Omega)\bigr).$ We say that $u$ is a weak solution of the equation $u_t = \dive(|\nabla u|^{p-2}\nabla u)$ in  $\Omega_T$ if
 \begin{equation*}
  \int\limits_0^T\!\!\!\int\limits_{\Omega}u\,\phi_t\,dxdt\,=\, \int\limits_0^T\!\!\!\int\limits_{\Omega}|\nabla u|^{p-2}\langle\nabla u,\nabla \phi\rangle\,dxdt\quad\text{when}\quad \phi \in C_0(\Omega_T).
 \end{equation*}
\end{definition}
\bigskip

Especially, $u\in L^2(\Omega_T)$ by the assumption.
The weak solutions for the regularized equation (\ref{regu})
are defined in a similar way, see
 (\ref{reg}). According to Theorem 4.2 on page 599 in [LSU] they have continuous second derivatives in all variables. We shall use the notation $\u$ for the solution of the regularized equation with boundary values $u$ on the parabolic boundary of $\Omega_T$. The boundary values are taken in the following sense:
 \begin{equation*}
\bullet\,\, \u-u \in L^p\bigl(0,T;W^{1,p}_0(\Omega)\bigr)\quad\text{and}\quad \bullet\,\,
 \lim_{\delta \to 0}\frac{1}{\delta} \int\limits_0^{\delta}\!\!\!\int\limits_{\Omega}|\u-u|^2\,dxdt\,=\,0.
 \end{equation*}

\section{The Time Derivative}

Our proof depends on the applicability of the rule
\begin{equation}\label{rule}
  \int\limits_0^T\!\!\!\int\limits_{\Omega}u\,\phi_t\,dxdt\,=\,- \int\limits_0^T\!\!\!\int\limits_{\Omega}\phi\,\nabla\!\cdot\!\bigl(|\nabla u|^{p-2}\nabla u\bigr)\,dxdt,
\end{equation}
when $\phi\in C_0^{\infty}(\Omega_T).$ Thus the theorem follows, if we can prove that the derivatives $\partial/ \partial x_j(|\nabla u|^{p-2}\nabla u)$ in the formula exist and belong to $L_{loc}^{2}(\Omega_T)$. Indeed, that we can do for $p >\tfrac{3}{2}$. Yet, for smaller values of $p,$ the negative exponent $p-2$ forces us to circumvent this expression, which is problematic when $\nabla u = 0.$ We use the regularized equation
\begin{equation}\label{reg}
  \int\limits_0^T\!\!\!\int\limits_{\Omega}\u\,\phi_t\,dxdt\,=\,- \int\limits_0^T\!\!\!\int\limits_{\Omega}\phi\,\nabla\!\cdot\!\Bigl(\!\bigl\{|\nabla \u|^2 +\e^2\bigr\}^{\frac{p-2}{2}}\nabla \u\Bigr)dxdt
    \end{equation}
    and prove that, as $\e \to 0,$  the derivatives
    $$\frac{\partial}{\partial x_{j}}\Bigl(\!\bigl\{|\nabla \u|^{2}+\e^2\bigr\}^{\!\frac{p-2}{2}}\frac{\partial\u}{\partial x_{k}}\Bigr)$$
    converge weakly in  $L_{loc}^{\theta}(\Omega_T)$ with some $\theta > 1$. Since $\u$ converges to $u$ locally in $L^2(\Omega_T)$ by Proposition \ref{conv}, the Theorem follows from the compactness result below, when we take into account that
    $$\left\vert \frac{\partial}{ \partial x_j}\bigl\{|\nabla \u|^2 +\e^2\bigr\}^{\frac{p-2}{2}}\frac{\partial \u}{\partial x_{k}}\right\vert \, \leq \,2 \bigl\{|\nabla \u|^2 +\e^2\bigr\}^{\frac{p-2}{2}}|\mathbb{D}^2\u|$$

      Assume that  $$u \in C\bigl(0,T;L^2(\Omega)\bigr) \cap L^p\bigl(0,T;W^{1,p}(\Omega)\bigr)$$ is a weak solution to $u_t = \dive(|\nabla u|^{p-2}\nabla u)$ in $\Omega_T$. Let $\u$ denote the solution of the regularized equation 
      $$  \frac{\partial \u}{\partial t}\,=\,\dive\left(\{|\nabla \u|^2 + \e^2\}^{\frac{p-2}{2}}\nabla \u\right)$$
      with the same boundary values as $u$ on the parabolic boundary of $\Omega_T$.

\bigskip
      \begin{lemma}\label{compact} We have uniformly with respect to $\e$:
        \begin{itemize}
        \item{$\mathbf{p\geq \frac{3}{2}.}$}\,\,
          \begin{equation*}
            \int\limits_0^T\!\!\!\int\limits_{\Omega}\zeta^2\bigl\{|\nabla \u|^2 +
            \e^2\bigr\}^{2(p-2)}|\mathbb{D}^2\u|^2dxdt\, \leq L\, <\,\infty,\qquad \e \leq 1. \qquad (\star)
          \end{equation*}
          \item{$\mathbf{1<p<\frac{3}{2}.}$}\,
            Under the extra assumption that $\|u\|_{\infty} < \infty$, the quantity
   $$ \int\limits_0^T\!\!\!\int\limits_{\Omega}\zeta^2\bigl\{|\nabla \u|^2 + \e^2\bigr\}^{\theta(p-2)}|\mathbb{D}^2\u|^2dxdt\, \leq L(\theta)\, <\,\infty,\qquad \e \leq 1.$$
                is uniformly bounded in $\e$ when $$1 < \theta < \frac{1}{2-p}.$$
               \item{$\mathbf{n=1.}$}\, In the one-dimensional case ($\star$) holds for all $p>1$.
        \end{itemize}
      \end{lemma}

      \bigskip
      \emph{Proof:} The second case is  Proposition $\ref{<1,5}$ and the two other cases are in Section \ref{Y}. Formally, ($\star$) is equation (2.16) in [Y].\quad $\Box$
        
\section{Convergence of the Approximation}
In this section we shun the extra assumption about the boundedness of the weak solution $u.$ This effort complicates the convergence proof for the $\u^{\,\,'}$s.
Recall the equations
\begin{equation*}
  \begin{cases}
    \frac{\partial u}{\partial t}\,\,=\,\dive\left(|\nabla u|^{p-2}\nabla u\right)\\
    \frac{\partial \u}{\partial t}\,=\,\dive\!\left(\!\{|\nabla \u|^2 + \e^2\}^{\frac{p-2}{2}}\nabla \u\right)

  \end{cases}
\end{equation*}
where $u = \u$ on the parabolic boundary of $\Omega_T.$

\bigskip

\begin{prop}\label{conv} Under the assumption
  $$u \in C(0,T;L^2(\Omega)) \cap L^p(0,T;W^{1,p}(\Omega))$$
  the convergence
  $$\u \to u \quad \text { in}\quad L^2(\Omega_T),\qquad \nabla \u \to \nabla u\quad \text { in}\quad L^p(\Omega_T)$$
  is valid.
\end{prop}

\bigskip

\emph{Proof:} 
Using the test function
$\phi = \u - u$ in both equations we get
\begin{align}\label{W}
&W_{\e} &\equiv& \int\limits_0^T\!\!\!\int\limits_{\Omega}\langle\{|\nabla \u|^2 + \e^2\}^{\frac{p-2}{2}}\nabla \u -|\nabla u|^{p-2}\nabla u,\nabla \u - \nabla u \rangle\,dxdt \nonumber\\
  &&=& - \frac{1}{2}\int\limits_{\Omega}(\u(x,T)-u(x,T))^2\,dx\,\leq 0.
\end{align}
Strictly speaking, in the equation for $u$ we must go via a time regularization; the Steklov average works well and the final inequality $W_{\e} \leq 0$ follows. Thus
\begin{align*}
  J_{\e}\, \equiv& \,  \int\limits_0^T\!\!\!\int\limits_{\Omega}\{|\nabla \u|^2 + \e^2\}^{\frac{p-2}{2}}|\nabla \u|^2\,dxdt\\
    \leq& - \int\limits_0^T\!\!\!\int\limits_{\Omega}|\nabla u|^p\,dxdt +  \int\limits_0^T\!\!\!\int\limits_{\Omega}| \nabla u|^{p-2}\langle \nabla u,\nabla \u \rangle \,dxdt\\
    &\quad+ \int\limits_0^T\!\!\!\int\limits_{\Omega}\{|\nabla \u|^2 + \e^2\}^{\frac{p-2}{2}}\langle \nabla \u,\nabla u \rangle \,dxdt.
\end{align*}
By Young's inequality the last integrand satisfies
\begin{equation*}
  \left|\{|\nabla \u|^2 + \e^2\}^{\frac{p-2}{2}}\langle \nabla \u,\nabla u\rangle \right| \leq \frac{1}{p}|\nabla u|^p + \frac{1}{q}\{|\nabla \u|^2 + \e^2\}^{\frac{p-2}{2}}|\nabla \u|^2,
\end{equation*}
where $q = p/(p-1).$ In the same way
$$ \left| \nabla u|^{p-2}\langle \nabla u,\nabla \u \rangle \right | \leq \frac{\sigma^p}{p}|\nabla \u|^p +\frac{1}{\sigma^q q}|\nabla u|^p,\qquad \sigma > 0.$$
Upon integration and absorption of a term, we arrive at
\begin{equation}\label{J1}
  J_{\e} \leq (p-1) (\sigma^{-q}-1)\int\limits_0^T\!\!\!\int\limits_{\Omega}|\nabla u|^p\,dxdt + \sigma^p\int\limits_0^T\!\!\!\int\limits_{\Omega}|\nabla \u|^p\,dxdt.
\end{equation}
In order to handle the last integral, we divide the domain  of integration into two parts: the set $|\nabla \u| \leq \e$ and the set $|\nabla \u| \geq \e.$ We have
\begin{align*}
  \int\limits_0^T\!\!\!\int\limits_{\Omega}|\nabla \u|^p\,dxdt\, &\leq \e^p\mathrm{mes}(\Omega_T) + \underset{|\nabla \u| \geq \e}{\int\!\!\int} |\nabla \u|^2|\nabla \u|^{p-2}dxdt\\
  &\leq \e^p\mathrm{mes}(\Omega_T) + 2^{\frac{2-p}{2}}\underset{|\nabla \u| \geq \e}{  \int\!\!\int}\{|\nabla \u|^2 + \e^2\}^{\frac{p-2}{2}}|\nabla \u|^2\,dxdt\\
  &\leq  \e^p\mathrm{mes}(\Omega_T) + 2^{\frac{2-p}{2}} J_{\e}.
\end{align*}
We insert this in equation (\ref{J1}) and obtain
\begin{equation*}
  J_{\e} \leq (p-1) (\sigma^{-q}-1)\int\limits_0^T\!\!\!\int\limits_{\Omega}|\nabla u|^p\,dxdt + \sigma^p \e^p\mathrm{mes}(\Omega_T) + \sigma^p 2^{\frac{2-p}{2}} J_{\e}.
\end{equation*}
We fix $\sigma > 0$ equal to a number, depending only on $p$, so small that the last $J_{\e}$-term can be absorbed into the left-hand side. It follows that
\begin{gather*}
  J_{\e} \leq C_p\biggl\{\int\limits_0^T\!\!\!\int\limits_{\Omega}|\nabla u|^p\,dxdt + \e^p\mathrm{mes}(\Omega_T) \biggr\},\\ \int\limits_0^T\!\!\!\int\limits_{\Omega}|\nabla  \u|^p\,dxdt \leq 3 C_p\biggl\{\int\limits_0^T\!\!\!\int\limits_{\Omega}|\nabla u|^p\,dxdt + \e^p\mathrm{mes}(\Omega_T) \biggr\}.
\end{gather*}
In particular, we have a uniform bound:
\begin{equation}\label{uni}
 \boxed{ \int\limits_0^T\!\!\!\int\limits_{\Omega}|\nabla  \u|^p\,dxdt \leq K, \qquad 0\leq \e \leq 1.}
  \end{equation}
Now we split $W_{\e}$ as
\begin{align*}
  W_{\e} &= \int\limits_0^T\!\!\!\int\limits_{\Omega}\langle|\nabla \u|^{p-2}\nabla  \u - |\nabla u|^{p-2}\nabla  u,\nabla \u - \nabla u \rangle\,dxdt\\
  &+ \int\limits_0^T\!\!\!\int\limits_{\Omega}\langle\left\{|\nabla \u|^2 +\e^2\right\}^{\frac{p-2}{2}}\nabla  \u - |\nabla u|^{p-2}\nabla  u,\nabla \u - \nabla u \rangle\,dxdt,
\end{align*}
and since $W_{\e} \leq 0$ by (\ref{W}),
\begin{align*}
  M_{\e} &\equiv \int\limits_0^T\!\!\!\int\limits_{\Omega}\langle|\nabla \u|^{p-2}\nabla  \u - |\nabla u|^{p-2}\nabla  u,\nabla \u - \nabla u \rangle\,dxdt\\
  &\leq \int\limits_0^T\!\!\!\int\limits_{\Omega}\langle|\nabla  \u|^{p-2}\nabla \u -  \left\{|\nabla \u|^2 +\e^2\right\}^{\frac{p-2}{2}}\nabla  \u,\nabla \u - \nabla u \rangle\,dxdt\, \equiv O_{\e}.
\end{align*}
We claim that $O_{\e} \rightarrow 0$ as $\e \to 0.$ Recall that  $1< p \leq 2$. Thus the inequality\footnote{\begin{gather*}0\leq|a|^{p-2}-(|a|^2+\e^2)^{\frac{p-2}{2}} = - \int_{0}^{1}\frac{\dif}{\dif t}(|a|^2+t\e^2)^{\frac{p-2}{2}}\,dt\\ = \frac{2-p}{2}\e^2 \int_{0}^{1}(|a|^2+t\e^2)^{\frac{p-4}{2}}\,dt \leq  \frac{2-p}{2}\e^2|a|^{p-4},\quad a\not = 0.\end{gather*}}
$$0\leq |a|^{p-2}-\bigl(|a|^2+\e^2\bigr)^{\frac{p-2}{2}} < \frac{2-p}{2}\e^2|a|^{p-2}\delta^{-2}\qquad |a| \geq \delta,$$
is available. Now we split the domain of integration for $O_{\e}$ into two parts and achieve
\begin{gather*}
  |O_{\e}| \leq \frac{2-p}{2}\e^2\delta^{-2}\underset{|\nabla \u| \geq \delta}{\int\!\!\int}|\nabla \u|^{p-1}|\nabla \u - \nabla u|\,dxdt\\
  +2\underset{|\nabla \u| \leq \delta}{\int\!\!\int}|\nabla \u|^{p-1}|\nabla \u - \nabla u|\,dxdt.
\end{gather*}
By H\"{o}lder's inequality
$$\int\!\!\!\int|\nabla \u|^{p-1}|\nabla \u - \nabla u|\,dxdt \leq \Bigl(\int\!\!\!\int|\nabla \u|^{p}\,dxdt\Bigr)^{\frac{p-1}{p}}\Bigl\{\|\nabla \u\|_{p}+ \|\nabla u\|_{p}\Bigr\}.$$
Recalling the uniform bound (\ref{uni}), we see that
\begin{gather*}
  |O_{\e}| \leq \tfrac{2-p}{2}\e^2\delta^{-2}K^{1-\frac{1}{p}}\!\left(K^{\frac{1}{p}}+\|\nabla u\|_p\right)+2\delta^{p-1}\! \left(K^{\frac{1}{p}}+\|\nabla u\|_p\right).
\end{gather*}
It follows that
$$\lim_{\e \rightarrow 0}O_{\e}\,=\,0.$$

The inequality\footnote{For vectors \begin{gather*}
\langle|b|^{p-2}b-|a|^{p-2}a,b-a\rangle\,\geq\,(p-1)|b-a|^2\bigl(1+|a|^2+|b|^2\bigr)^{\frac{p-2}{2}},\quad 1<p\leq2.
\end{gather*} }
\begin{equation}
  \int\limits_0^T\!\!\!\int\limits_{\Omega}\frac{|\nabla \u - \nabla u|^2\,dxdt}{\bigl(1+ |\nabla u|^2+|\nabla \u|^2\bigr)^{\frac{2-p}{2}}}\, \leq \, M_{\e}\, \leq \, O_{\e}
\end{equation}
shows in combination with
\begin{gather*}
  \int\limits_0^T\!\!\!\int\limits_{\Omega}|\nabla \u - \nabla u|^p\,dxdt =   \int\limits_0^T\!\!\!\int\limits_{\Omega}|\nabla \u - \nabla u|^p\biggl\{\frac{1+|\nabla u|^2+|\nabla \u|^2}{1+|\nabla u|^2+|\nabla \u|^2}\biggr\}^{\!\tfrac{p(2-p)}{4}}dxdt\\
  \leq \biggl\{  \int\limits_0^T\!\!\!\int\limits_{\Omega}\frac{|\nabla \u - \nabla u|^2\,dxdt}{\bigl(1+|\nabla u|^2+|\nabla \u|^2\bigr)^{\frac{2-p}{2}}}\biggr\}^{\!\tfrac{p}{2}} \int\limits_0^T\!\!\!\int\limits_{\Omega}\bigl(1+|\nabla u|^2+|\nabla \u|^2\bigr)^{\frac{p}{2}}\,dxdt
\end{gather*}
and the uniform bound (\ref{uni})
that
\begin{equation}
  \lim_{\e \to 0} \int\limits_0^T\!\!\!\int\limits_{\Omega}|\nabla \u - \nabla u|^p\,dxdt =  0.
\end{equation}

The convergence $\u \to u$ in $L^2(\Omega_T)$ can be extracted from the above proof,
according to which
$$\frac{1}{2}\int\limits_{\Omega}\bigl(\u(x,T)-u(x,T)\bigr)^2\,dx = - W_{\e} = O_{\e}-M_{\e}\leq O_{\e} \to 0.$$
When we replace $T$ by $t,\,\,0 < t < T$, the same bound as before will majorize $O_{\e}$ \emph{simultaneously} for all $t$. Integrating with respect to $t,$ we obtain
$$\int\limits_0^T\!\!\!\int\limits_{\Omega}\bigl(\u(x,t)-u(x,t)\bigr)^2\,dxdt\, \leq \,2TO_{\e}.$$
This concludes the convergence proof.\qquad $\Box$

\section{The Main Identity}\label{Y}

In order to derive estimates for the derivatives
$$\frac{\partial}{\partial x_{j}}\bigl(|\nabla u|^{p-2}\nabla u \bigr)$$
we differentate the \emph{regularized} equation
$$\frac{\partial \u}{\partial t} = \dive\! \left(\!\bigl(|\nabla \u|^2 + \e^2\bigr)^{\frac{p-2}{2}}\nabla \u\!\right).$$
Using the abbreviations
\begin{gather*}
  u_{\e,j} = \frac{\partial}{\partial_{x_j}}\u,\quad \v = |\nabla \u|^2, \quad
\V = |\nabla \u|^2 + \e^2 = \v^2 + \e^2
\end{gather*}
we have
\begin{equation}\label{differentiated}
  \frac{\partial}{\partial t}u_{\e,j} = \dive\!\left(\V^{\frac{p-2}{2}}\nabla u_{\e,j} + \nabla u_{\e,j} \frac{\partial}{\partial_{x_j}}\V^{\frac{p-2}{2}}\right).
  \end{equation}
We note
\begin{gather*}
  \frac{\partial}{\partial_{x_j}}\v = 2\langle \nabla \u,\nabla u_{\e,j}\rangle    ,\qquad |\nabla \v|^2 \leq 4|\nabla \u|^2|\mathbb{D}^2\u|^2 \\
  \frac{\partial}{\partial_{x_j}}\V^{\frac{p-2}{2}} = (p-2)\V^{\frac{p-4}{2}}\langle \nabla \u,\nabla u_{\e,j}\rangle.
\end{gather*}
In weak form the equation becomes
\begin{align}\label{start}
 - &\int\limits_0^T\!\!\!\int\limits_{\Omega}\phi_j\frac{\partial u_{\e,j}}{\partial t} \, dxdt\\ =
  &\int\limits_0^T\!\!\!\int\limits_{\Omega}\left(\V^{\frac{p-2}{2}}\langle\nabla u_{\e,j},\nabla \phi_j\rangle +  (p-2)\V^{\frac{p-4}{2}}\langle \nabla \u,\nabla u_{\e,j}\rangle \langle \nabla \u, \nabla \phi_j\rangle\!\right)dxdt,\nonumber
\end{align}
valid at least for all test functions $\phi_j \in C_0^{\infty}(\Omega_T),\,\, j= 1,2,\dots,n.$ (In fact, it is not needed that $\phi_j =0$ when $t=0$ or $t=T.$) We use the test functions
$$\phi_j= \zeta^2\V^{\alpha}u_{\e,j},\qquad \zeta\in C_0^{\infty}(\Omega_T)$$
and sum the formulas to reach the identity below. (Such identities often serve to derive Caccioppoli inequalities.) We shall keep $1-p<2\alpha<0.$ Always, $0\leq \zeta \leq 1.$
\paragraph{Fundamental formula}
\begin{align} &\int\limits_0^T\!\!\!\int\limits_{\Omega}\zeta^2\V^{\frac{p-2+2\alpha}{2}}|\mathbb{D}^2\u|^2\,dxdt\qquad\qquad\text{\textsf{Main Term}} \tag{I}\\
  +&\,\frac{p-2+2\alpha}{4} \int\limits_0^T\!\!\!\int\limits_{\Omega}\zeta^2\V^{\frac{p-2+2\alpha}{2}-1}|\nabla \v|^2\,dxdt \tag{I\!I}\\
  +&\,\frac{\alpha (p-2)}{2}\int\limits_0^T\!\!\!\int\limits_{\Omega}\zeta^2\V^{\frac{p-2+2\alpha}{2}-2}\langle \nabla \u,\nabla \v \rangle^2\,dxdt \tag{I\!I\!I}\\
  +&\,\frac{1}{2(\alpha +1)}\bigg\lbrack\int\limits_{\Omega}\zeta^2\V^{\alpha+1}\,dx\bigg\rbrack_{0}^{T} \tag{I\!V}\\
  =&\,\, (2-p)\int\limits_0^T\!\!\!\int\limits_{\Omega}\zeta\V^{\frac{p-2+2\alpha}{2}-1}\langle \nabla \u,\nabla \v \rangle\langle \nabla \zeta,\nabla \u\rangle\,dxdt \tag{V}\\
  - &\int\limits_0^T\!\!\!\int\limits_{\Omega}\zeta\V^{\frac{p-2+2\alpha}{2}}\langle \nabla \zeta, \nabla \u \rangle\,dxdt \tag{V\!I}\\
  +&\,\frac{1}{\alpha+1}\int\limits_0^T\!\!\!\int\limits_{\Omega}\V^{\alpha +1}\zeta \zeta_{t}\,dxdt \tag{V\!I\!I}
\end{align}

\bigskip
The proof is a straightforward calculation. (Compare with formula (2.5) in [Y] and formula (2.20) on page 166 in [WZY].) We only mention how to treat the part with the time derivative:
\begin{align*}
  &\phi_j\,\frac{\partial}{\partial t}u_{\e,j} = \zeta^2\V^{\alpha}\,\frac{\partial}{\partial t}\Bigl(\frac{u_{\e,j}^2}{2}\Bigr)\\
  &\zeta^2\V^{\alpha}\,\frac{\partial}{\partial t}\Bigl(\frac{\v}{2}\Bigr) = \frac{1}{2}\zeta^2\,\frac{\partial}{\partial t}\Bigl(\frac{\V^{\alpha+1}}{\alpha +1}\Bigr).
\end{align*}
Thus, upon summation, the left-hand side of (\ref{start}) becomes
$$ -\sum_{j=1}^{n}\int\limits_0^T\!\!\!\int\limits_{\Omega}\phi_j\frac{\partial u_{\e,j}}{\partial t}\,dxdt = \frac{1}{2}\bigg\lbrack\int\limits_{\Omega}\zeta^2\frac{\V^{\alpha+1}}{\alpha +1}\,dx\bigg\rbrack_0^T- \int\limits_0^T\!\!\!\int\limits_{\Omega}\zeta\zeta_t\frac{\V^{\alpha+1}}{\alpha +1}\,dxdt.$$
The right-hand side yields six terms, since the right-hand side of (\ref{differentiated}) is multiplied by
$$\nabla \phi_j = \zeta^2\V^{\alpha}\nabla u_{\e.j}+ \alpha\zeta^2\v^{\alpha-1}u_{\e,j}\nabla\V+2\V^{\alpha}u_{\e,j}\zeta\nabla\zeta;$$
 two  similar terms are joined in  term I\!I.
    
 Always, $0 > 2\alpha > 1-p$ and $1 < p \leq 2$, which means that \emph{the factor in front of term} I\!I\emph{ is negative}. The integral itself is of the same magnitude as  term I, and
 \begin{equation}\label{2}
   |\nabla \v|^2 \leq 4 \v |\mathbb{D}^2\u|^2\leq 4\V |\mathbb{D}^2\u|^2.
   \end{equation}
 This causes the \emph{constraint}: $p-1+2\alpha > 0.$
Term I\!I\!I is positive, but since the expression
$$\langle \nabla \u, \nabla \v \rangle^2 = 4\sum_{i,j=1}^{n}\frac{\partial \u}{\partial x_i}
\frac{\partial \u}{\partial x_j}\frac{\partial^2 \u}{\partial x_i \partial x_j}$$
  may vanish, it is of little use, \emph{except in the one dimensional case} when term I\!I\!I matches term I. 
  
  \subsection*{Estimation of some terms}
 
  \paragraph{Vanishing of term I\!V} It is zero, as $\zeta$ has compact support also in the time direction.
  \paragraph{Absorption of term V} We can use Young's inequality to absorb term V into the main term I. Now by (\ref{2})
  $$|\langle\nabla \u,\nabla \v\rangle\langle\nabla \zeta,\nabla \u\rangle| \leq |\nabla \u|^2|\nabla \zeta|\,2\,\V^{\frac{1}{2}}|\mathbb{D}\u|$$
  and with a small parameter $\sigma > 0$
  \begin{equation}
    |\mathrm{V}| \leq (2-p)\sigma \mathrm{I} + (2-p)\sigma^{-1}\int\limits_0^T\!\!\!\int\limits_{\Omega}\V^{\frac{p+2\alpha}{2}}|\nabla \zeta|^2\,dxdt.
  \end{equation}
  \paragraph{Term V\!I}
  Since $|\langle \nabla \zeta, \nabla \v\rangle| \leq 2|\nabla \zeta||\nabla \u||\mathbb{D}^2\u|$, we get the same as above:
  \begin{equation}
    |\mathrm{V\!I}| \leq \sigma \mathrm{I} + \sigma^{-1}\int\limits_0^T\!\!\!\int\limits_{\Omega}\V^{\frac{p+2\alpha}{2}}|\nabla \zeta|^2\,dxdt.
  \end{equation}

  With these arrangements the main formula yields the estimate 
  \begin{align}
    &(1-(3-p)\sigma)\mathrm{I}\, +\, \mathrm{I\!I}\, +\, \mathrm{I\!I\!I} \\
    &\leq (3-p) \sigma^{-1}\int\limits_0^T\!\!\!\int\limits_{\Omega}\V^{\frac{p+2\alpha}{2}}|\nabla \zeta|^2\,dxdt + \mathrm{V\!I\!I}\nonumber
  \end{align}

  \subsection*{The one-dimensional case}

  In one space dimension we have
  $$\u' = \frac{\partial \u}{\partial x},\quad\u'' = \frac{\partial^2 \u}{\partial x ^2}  \quad \v = \u'^{\,2}, \quad \v' = 2\u'\u''.$$
  We fix $2\alpha = p-2,$ which is negative. Then the sum $\mathrm{I} + \mathrm{I\!I} + \mathrm{I\!I\!I}$ can be written as
  \begin{equation*} 
    \int\limits_0^T\!\!\!\int\limits_{\Omega}\zeta^2\V^{p-2}\bigl|\frac{\partial^2 \u}{\partial x^2}\bigr|^2 \left\{1+2(p-1)\u'^{\,2}\v^{-1}+(p-2)^2\u'^{\,4}\v^{-2}\right\}dxdt.
  \end{equation*}
  The expression in braces is a perfect square and can be estimated as
  $$\{1+\dots \v^{-2}\} = \left(\frac{(p-1)\u'^{\,2}+\e^2}{\u'^{\,2}+\e^2}\right)^2 \geq (p-1)^2.$$
  Thus the total estimate in one dimension reads
  \begin{align*}
    \bigl((p-1)^2 -(3-p)\sigma\bigr) \int\limits_0^T\!\!\!\int\limits_{\Omega}\zeta^2\V^{p-2}\bigl|\frac{\partial^2 \u}{\partial x^2}\bigr|^2dxdt\\
    \leq (3-p) \sigma^{-1}\int\limits_0^T\!\!\!\int\limits_{\Omega}\V^{p-1}|\nabla \zeta|^2\,dxdt + \frac{2}{p}\int\limits_0^T\!\!\!\int\limits_{\Omega}\V^{\frac{p}{2}}|\zeta\zeta_t|\,dxdt
  \end{align*}
  Now we only have to fix $\sigma$ small enough, noticing that
  $$\V^{p-1} \leq \u^{p} +1, \qquad \V^{\frac{p}{2}} \leq 2(\u^{p}+\e^p),$$
  to obtain the majorant
  \begin{equation}\label{onedim}
    \int\limits_0^T\!\!\!\int\limits_{\Omega}\zeta^2\V^{p-2}\bigl|\frac{\partial^2 \u}{\partial x^2}\bigr|^2dxdt\leq C(p) \biggl\{ \int\limits_0^T\!\!\!\int\limits_{\Omega}|\nabla \u|^p\,dxdt + 1 \biggr\}.
  \end{equation}
  The majorant is finite and, by (\ref{uni}) independent of $\e$, but the constant factor $C(p)$ depends also on $\|\zeta_t\|_{\infty}.$

  To proceed, use
  $$\left\vert\frac{\partial}{\partial x}\bigl\{\u'^{2}+\e^2\bigr\}^{\frac{p-2}{2}}\u'\right\vert\,\leq\, (p-1)\V^{\frac{p-2}{2}}\bigl|\frac{\partial^2 \u}{\partial x^2}\bigr|^2$$
  to conclude that
  \begin{equation*}
    \frac{\partial\,}{\partial x}\bigl\{\u'^{2}+\e^2\bigr\}^{\frac{p-2}{2}}\u'\quad\text{converges weakly in}\quad L^2_{loc}(\Omega_T)
  \end{equation*}
  at least through a subsequence. Thus we may pass to the limit under the integral signs in
  $$-\int\limits_0^T\!\!\!\int\limits_{\Omega}\u\frac{\partial\phi}{\partial t}\,dxdt = \int\limits_0^T\!\!\!\int\limits_{\Omega}\phi\, \frac{\partial}{\partial x}\Bigl(\!\{\u'^{2}+\e^2\}^{\frac{p-2}{2}}\u'\Bigr)dxdt$$
  and conclude that \emph{the time derivative} $u_t$ \emph{exists and belongs locally to} $L^2$. The limit is some function.

  \subsection*{General Estimate,\,\,$1<p<2$}

  In several space dimensions term I\!I\!I is no longer so useful, so one may as well skip it since it is positive when $\alpha <0.$  However, it is convenient to use it to counterbalance a portion of term V:
  $$ |\mathrm{V}| \leq \mathrm{I\!I\!I} + \frac{2-p}{|\alpha|} \int\limits_0^T\!\!\!\int\limits_{\Omega}\V^{\frac{p+2\alpha}{2}}
  |\nabla \zeta|^2\,dxdt,$$
  where Young's inequality was used. Now we have the general estimate
  \bigskip
  \begin{lemma}[$1<p<2.$] Let $\sigma > 0.$  We have
  \begin{align}\label{gen}
    &\bigl(p-1+2\alpha -\sigma)\int\limits_0^T\!\!\!\int\limits_{\Omega}\zeta^2\V^{\frac{p-2+2\alpha}{2}}|\mathbb{D}^2\u|^2\,dxdt \\
    &\leq\, \Bigl(\sigma^{-1} +\frac{2-p}{|\alpha|}\Bigr) \int\limits_0^T\!\!\!\int\limits_{\Omega}\V^{\frac{p+2\alpha}{2}}
    |\nabla \zeta|^2\,dxdt + \frac{1}{\alpha+1}\int\limits_0^T\!\!\!\int\limits_{\Omega}\V^{\alpha+1}\zeta\zeta_t\,dxdt.\nonumber
  \end{align}
  \end{lemma}
  \bigskip
  This is worthless if one does not obey the \begin{equation*} \boxed{\text{ restriction:}\qquad
      p-1+2\alpha\, >\, 0}\end{equation*}
  
 \section{The case $\frac{3}{2} < p <2$}
  
  Again, we take  $2\alpha = p-2$. Then by (\ref{gen})
   \begin{align*}
    \bigl(2p-3& -\sigma\bigr)\int\limits_0^T\!\!\!\int\limits_{\Omega}\zeta^2\V^{p-2}|\mathbb{D}^2\u|^2\,dxdt \\
    \leq& \left(\sigma^{-1} + 2\right) \int\limits_0^T\!\!\!\int\limits_{\Omega}\V^{p-1}
    |\nabla \zeta|^2\,dxdt + \frac{1}{p-1}\int\limits_0^T\!\!\!\int\limits_{\Omega}\V^{\frac{p}{2}}\zeta\zeta_t\,dxdt,
   \end{align*}
   Provided that $2p > 3$, this yields the desired local bound with a majorant free of $\e$ according to (\ref{uni}). 
   
   \section{An ''Energy Term'' with $p<\frac{3}{2}$}

   In the demanding case $p < \tfrac{3}{2}$ we need to estimate  the last integral in (\ref{gen}). In this case,\emph{ we assume that the solution is bounded}:$\, \|u\|_{\infty} <\infty.$
   Obviously
   $$\V^{\alpha+1}= \V^{\alpha}\bigl(\e^2+|\nabla \u|^2\bigr) \leq \e^{2(\alpha+1)} + \V^{\alpha}\langle\nabla\u,\nabla\u\rangle,$$
   since $\alpha < 0.$ If $\zeta \in C^{\infty}_{0}(\Omega)$, then\footnote{It is essential that this holds also when $\zeta\zeta_t <0$.}
   $$\int\limits_0^T\!\!\!\int\limits_{\Omega}\zeta\zeta_t\V^{\alpha+1}\,dxdt\leq \e^{2(\alpha+1)}\!\int\limits_0^T\!\!\!\int\limits_{\Omega}|\zeta\zeta_t|\,dxdt + \int\limits_0^T\!\!\!\int\limits_{\Omega}\zeta\zeta_t\V^{\alpha}\langle\nabla\u,\nabla\u\rangle\,dxdt.$$
   Integration by parts yields
   \begin{align*}     \int\limits_0^T\!\!\!\int\limits_{\Omega}\zeta\zeta_t\V^{\alpha}\langle\nabla\u,\nabla\u\rangle\,dxdt = -\, \int\limits_0^T\!\!\!\int\limits_{\Omega}\u\nabla\!\cdot\!\bigl(\zeta\zeta_t\V^{\alpha}\nabla \u\bigr)dxdt\\ =
     - \! \int\limits_0^T\!\!\!\int\limits_{\Omega}\!\left(\u\V^{\alpha}\langle\nabla( \zeta\zeta_t),\!\nabla \u\rangle - 2\alpha\zeta\zeta_t\u\V^{\alpha-1}\!\langle\nabla \u,\!\nabla u_{\e,j}\rangle - \u\zeta\zeta_t\V^{\alpha}\Delta \u\right)\!dxdt\\
     \leq \int\limits_0^T\!\!\!\int\limits_{\Omega}|\u\V^{\alpha}\nabla(\zeta\zeta_t)||\nabla \u|\,dxdt +(1+2|\alpha|)\int\limits_0^T\!\!\!\int\limits_{\Omega}|\zeta\zeta_t\u|\V^{\alpha}|\mathbb{D}^2\u|\,dxdt.
   \end{align*}
   The last integral can now be absorbed into the main term in (\ref{gen}). To see this,
   factorize
   $$\zeta\zeta_t\V^{\alpha} = \zeta\V^{\frac{p+2\alpha-2}{4}}\!\cdot\zeta_t\V^{\frac{2-p + 2\alpha}{4}}$$
   and select a small $\kappa >0$ for Young's inequality $2ab\leq \kappa a^2+\kappa^{-1}b^2.$ 
   We arrive at the final bound below, after some arrangements.

   \bigskip
   \begin{lemma}[Energy Estimate]\label{energy}
     We have
     \begin{align*} \overbrace{\int\limits_0^T\!\!\!\int\limits_{\Omega}\!\zeta\zeta_t\V^{\alpha+1}dxdt}^{(\alpha+1)\times\text{Term VII}}\leq \e^{2(\alpha+1)}\!\int\limits_0^T\!\!\!\int\limits_{\Omega}\!|\zeta\zeta_t|dxdt +2\|\u\|_{\infty} \! \int\limits_0^T\!\!\!\int\limits_{\Omega}|\nabla (\zeta\zeta_t)|\V^{\frac{2\alpha+1}{2}}dxdt\nonumber\\
       \frac{1+2|\alpha|}{2}\biggl\{\!\kappa\underbrace{\!\int\limits_0^T\!\!\!\int\limits_{\Omega}\zeta^2\V^{\frac{p+2\alpha-2}{2}}|\mathbb{D}^2\u|^2\,dxdt}_{\text{Term I}}
       +\kappa^{-1}\int\limits_0^T\!\!\!\int\limits_{\Omega}\zeta_t^2\V^{\frac{2-p+2\alpha}{2}}\,dxdt\!\biggr\}\|\u\|_{\infty}\nonumber
     \end{align*}
\end{lemma}
     \section{The Case $p < \tfrac{3}{2}$}

     We combine the estimate in Lemma \ref{energy} with the general inequality (\ref{gen})
     writing
     $$p-2 +2\alpha = \theta(p-2),\quad 2\alpha = (\theta-1)(p-2) > 1-p$$
     We must obey the restriction $p-1+2\alpha > 0$, which means that
     \begin{equation*}
       \boxed{1<\theta<\frac{1}{2-p}}
     \end{equation*}
     We obtain 
     \begin{align*}
       \Bigl(1+\theta(p-2)-&\sigma -\frac{2\kappa\|\u\|_{\infty}}{2-p}\Bigr)\int\limits_0^T\!\!\!\int\limits_{\Omega}\zeta^2\V^{\theta(p-2)}|\mathbb{D}^2\u|^2\,dxdt\\
       \leq \frac{\e^{2(\alpha+1)}}{\alpha+1}&\int\limits_0^T\!\!\!\int\limits_{\Omega}|\zeta\zeta_t|\,dxdt
         +\frac{2\|\u\|_{\infty}}{\alpha+1}\int\limits_0^T\!\!\!\int\limits_{\Omega}|\nabla(\zeta\zeta_t)|\V^{\frac{2\alpha+1}{2}}dxdt\\
           +& \frac{2\kappa^{-1}\|\u\|_{\infty}}{2-p}\int\limits_0^T\!\!\!\int\limits_{\Omega}\zeta^2_t\V^{\frac{2-p+2\alpha}{2}}dxdt\\
           +&\Bigl(\sigma^{-1}+\frac{2-p}{|\alpha|}\Bigr)\int\limits_0^T\!\!\!\int\limits_{\Omega}\V^{\frac{p+2\alpha}{2}}|\nabla \zeta|^2\,dxdt,
     \end{align*}
     where simplifying estimations like
    $\frac{1+2|\alpha|}{2(\alpha+1)} \leq \frac{2}{2-p}$
     have been made. Notice that $\alpha+1 > 0.$ Now the powers of $\V$ are decisive; they must be positive and no greater than $p/2.$ Our permanent restriction $0>2\alpha>1-p$ leads to
     \begin{gather*}0 < \tfrac{2\alpha+1}{2} < \tfrac{p}{2},\quad\tfrac{1}{2} < \tfrac{p +2\alpha}{2} <  \tfrac{p}{2},\quad
       \tfrac{3-2p}{2}< \tfrac{2-p+2\alpha}{2} <  \tfrac{p}{2},
     \end{gather*}
     but now we need $p< \tfrac{3}{2}$ in order to assure 
     $$\frac{2-p+2\alpha}{2} > 0.$$
     We see that the exponents are in the right range, and for the last three integrals we can use
     $$\V^{\beta}\leq |\nabla \u|^p + 1 \quad\text{for}\quad 0 < \beta < \tfrac{p}{2},\quad 0<\e\leq 1$$ and then the uniform bound (\ref{uni}). 
     By the Maximum Principle $\|\u\|_{\infty}\leq \|u\|_{\infty}$ (equality holds). Hence we have the following result:

     \bigskip

     \begin{prop} \label{<1,5} Let $1 < p <\frac{3}{3}.$ Fix $\theta$ in the range $1 < \theta < \frac{1}{2-p}.$ Then the integral
       $$\boxed{\int\limits_0^T\!\!\!\int\limits_{\Omega}\zeta^2\V^{\theta(p-2)}|\mathbb{D}^2\u|^2\,dxdt\, \leq\, L(\theta), \quad 0<\e\leq1,}$$
 is uniformly bounded in $\e$.  The bound $L(\theta)$ depends also on $\, p,\,\|u\|_{\infty},\,\| \nabla \zeta\|_{\infty},$\\$\|\zeta_t\|_{\infty},$ and the constant $K$ in (\ref{uni}).
       \end{prop}

\end{document}